\documentstyle[12pt,amssymb,a4]{article}

\def\mbo#1{{\mathchoice {\mbox{\normalshape\normalsize #1}}
{\mbox{\normalshape\normalsize
#1}} {\mbox{\normalshape\scriptsize #1}} {\mbox{\normalshape\tiny #1}} }}

\def\omi#1{\buildrel { #1 \atop ^\vee} \over .}

\def\bbbone{{\mathchoice {\mbox{$1\mskip-4mu$l}} {\mbox{$1\mskip-4mu$l}}
{\mbox{$1\mskip-4.5mu$l}} {\mbox{$1\mskip-5mu$l}} }}
\def\gR{{\Bbb R}}

\def\gC{{\Bbb C}}
\def\gZ{{\Bbb Z}}
\def\gN{{\Bbb N}}

\def\cA{{\cal A}}
\def\cB{{\cal B}}
\def\cC{{\cal C}}
\def\cF{{\cal F}}
\def\cG{{\cal G}}
\def\cH{{\cal H}}
\def\cI{{\cal I}}

\def\cM{{\cal M}}
\def\cN{{\cal N}}

\def\cQ{{\cal Q}}
\def\cS{{\cal S}}

\def\cZ{{\cal Z}}
\def\fg{{\frak g}}
\def\fh{{\frak h}}
\def\Im{{\mbo{Im}}}
\def\ker{{\mbo{Ker}}}
\def\der{{\mbo{Der} }}
\def\out{{\mbo{Out} }}
\def\Int{{\mbo{Int} }}
\def\End{{\mbo{End} }}

\newtheorem{proposition}{Proposition}[section]
\newtheorem{definition}{Definition}[section]
\newtheorem{lemme}{Lemma}[section]

\begin{document}

\vspace*{3cm}
\begin{center}
{\large\bf SUBMANIFOLDS AND QUOTIENT MANIFOLDS\\
\bigskip
IN NONCOMMUTATIVE GEOMETRY}
\end{center}
\vspace{1cm}

\begin{center}
Thierry MASSON\\
\vspace{0.3cm}
{\small Laboratoire de Physique Th\'eorique et Hautes
Energies\footnote{Laboratoire associ\'e au Centre National de la
Recherche Scientifique - URA D0063}\\
Universit\'e Paris XI, B\^atiment 211\\
91 405 Orsay Cedex, France\\
e-mail: masson@qcd.th.u-psud.fr}
\end{center}
\vspace{1cm}

\begin{center}
\today
\end{center}

\vspace{1cm}
\begin{abstract}
We define and study noncommutative generalizations of submanifolds and
quotient manifolds, for the derivation-based differential calculus
introduced by M.~Dubois-Violette and P.~Michor. We give examples to
illustrate these definitions.
\end{abstract}

\vspace {1cm}

\noindent L.P.T.H.E.-ORSAY 95/53

\newpage

Various noncommutative generalizations of differential forms have been
proposed as well as generalizations of vector bundles and connections.
What is still missing is the concept of a submanifold and of a quotient
manifold, that is, how the differential structure of a given algebra
must be related to the differential structure of a subalgebra
(``quotient manifold'') or a quotient algebra (``submanifold''). In this
paper, we propose a definition of a noncommutative submanifold and of a
noncommutative quotient manifold within the context of the
derivation-based differential calculus first introduced by
M.~Dubois-Violette \cite{MDV1}, and completed
\cite{DVM1,DVM2} with P.~Michor.

In the first Section, we recall various definitions related to this
differential calculus.  In the second Section, we recall the definition
of Hochschild cohomology and other cohomologies which will be used
later. Submanifolds and quotient manifolds are defined respectively in
Sections 3 and 4.

\section{Noncommutative Differential Structures}

In noncommutative geometry, the algebra of smooth functions on a
manifold is replaced by a noncommutative algebra (See, for example
\cite{CO1}, \cite{MDV2}). Geometric objects are first expressed in terms
of the algebra of functions and then they can be generalized to the
noncommutative case. In this section, we recall the definition of
differential forms, central bimodules and connections as they are given
in \cite{MDV1},\cite{DVM1} and \cite{DVM2}.

\subsection{Noncommutative Differential Forms}

Let $\cA$ denote an associative algebra with unit. It is then the
generalization of the algebra of smooth functions on a compact manifold.
The center of the algebra will be denoted by $\cZ(\cA)$. The
differential forms we wish to introduce are based on derivations, the
algebraic generalizations of vector fields on a manifold:
\[ \der(\cA) = \{ X: \cA \rightarrow \cA \,\, /\,\,  X(ab) = X(a)b + aX(b) \}
\]
$\der(\cA)$ is naturally a $\cZ(\cA)$-module and a Lie algebra.

The two noncommutative generalizations of the graded differential
algebra of differential forms which we shall need \cite{MDV1,DVM1} are
constructed as follows. Let $C_{\cZ(\cA)}(\der(\cA) , \cA)$ be the
graded algebra of antisymmetric $\cZ(\cA)$-multi\-linear mappings from
$\der(\cA)$ to $\cA$. Notice that this algebra is not graded
commutative. In degree 0 we take $C_{\cZ(\cA)}^0(\der(\cA) , \cA) =
\cA$. We introduce a differential $d$ by the Koszul formula:
\begin{eqnarray*}
d\omega(X_1, \dots , X_{n+1}) &=& \sum_{i=1}^{n+1} (-1)^{i+1} X_i
\omega( X_1, \dots \omi{i} \dots, X_{n+1}) \nonumber\\ & & + \sum_{1\leq
i < j \leq n+1} (-1)^{i+j} \omega( [X_i, X_j], \dots \omi{i} \dots
\omi{j} \dots , X_{n+1})
\end{eqnarray*}
for any $\omega \in C_{\cZ(\cA)}^n(\der(\cA) , \cA)$ and any set of
derivations $X_i$.

\bigskip
Now, we can introduce the first generalization of differential forms
over $\cA$. We define $\Omega_{\der}(\cA)$ to be the smallest
differential graded subalgebra of the algebra $C_{\cZ(\cA)}(\der(\cA) ,
\cA)$ which contains $\cA$. Every element $\omega\in
\Omega_{\der}^n(\cA)$ can be written as a finite sum of elements of the
type $a_0 da_1 \dots da_n$, where $da\in \Omega_{\der}^1(\cA)$ is the
1-form $X\mapsto Xa \in\cA$, and where the product is that of
$C_{\cZ(\cA)}(\der(\cA) , \cA)$.

The second differential graded algebra of forms we shall use is the
algebra $C_{\cZ(\cA)}(\der(\cA) , \cA)$ itself, denoted by
$\underline{\Omega}_{\der}(\cA)$. We refer the reader to \cite{DVM2} for
the relationship between $\Omega_{\der}(\cA)$, $\der(\cA)$ and
$\underline{\Omega}_{\der}(\cA)$ from the point of view of duality .

\bigskip
There is a canonical Cartan operation $i_X$ of the Lie algebra
$\der(\cA)$ on $\Omega_{\der}(\cA)$ and $\underline{\Omega}_{\der}(\cA)$
\cite{MDV1}. For any $X\in \der(\cA)$, one defines the antiderivation of
degree $-1$
\[ i_X : \Omega_{\der}^n(\cA) \rightarrow \Omega_{\der}^{n-1}(\cA) \]
by
\[ (i_X\omega)(X_1, \dots, X_{n-1}) = \omega(X,X_1, \dots, X_{n-1}) \]
with $i_X a = 0$ for any $a\in\cA = \Omega_{\der}^0(\cA)$. It follows
that $L_X = i_X d + di_X$ is a derivation of degree 0 on
$\Omega_{\der}(\cA)$.

\subsection{Central Bimodules and Connections}

In ordinary differential geometry, vector bundles of finite rank can be
considered from an algebraic point of view through their space of
sections. In fact this space is a finite projective module over the
algebra of smooth functions. In noncommutative geometry, the
generalization of a vector bundle will then be such a module over the
algebra. But, since $\cA$ is noncommutative this can be a right module,
a left module or a bimodule.

In \cite{DVM1}, it was proposed that this generalization should at least
have the structure of a central bimodule. We recall that a central
bimodule $\cM$ is a bimodule over $\cA$ and which is also a module over
the center $\cZ(\cA)$ of $\cA$ in the commutative sense. That is, for
any $z\in \cZ(\cA)$ and $m\in\cM$, one has $zm=mz$.

\bigskip
It is then easy to introduce the notion of a connection on a central
bimodule. A connection on $\cM$ is a linear mapping $\nabla$ from
$\der(\cA)$ into the linear endomorphisms of $\cM$ such that
\begin{eqnarray*}
\nabla_{zX} m &=& z\nabla_X m \\
\nabla_X (amb) &=& (Xa)mb + a(\nabla_X m )b + am(Xb)
\end{eqnarray*}
for any $X\in\der(\cA)$, $z\in \cZ(\cA)$, $a,b\in\cA$ and $m\in\cM$.

The curvature of this connection is defined by the usual formula
\[ R(X,Y) = \nabla_X\nabla_Y - \nabla_Y\nabla_X - \nabla_{[X,Y]} \]
for any $X,Y\in\der(\cA)$. $R(X,Y)$ is an $\cA$-bimodule endomorphism of
$\cM$, antisymmetric and $\cZ(\cA)$-linear in $X,Y$. We refer the reader
to \cite{DVM2} for more properties on these connections.

\section{Hochschild Cohomology and Related Cohomologies}

In this section we introduce a class of subcomplexes of the Hochschild
complex of an associative algebra, and their cohomology. These
cohomologies, in degree 1, will be useful in the next section.

\subsection{Hochschild Cohomology}

We recall the definition of the ordinary Hochschild cohomology. Let
$\cA$ be an associative algebra with unit over $\gC$, and $\cM$ a
bimodule over $\cA$.

We define the complex $C(\cA;\cM)$ as follows: $C^n(\cA;\cM)$ is the
linear space of $\gC$-linear mappings from $\cA^{\otimes^n}$ to $\cM$.
In degree 0, we set $C^0(\cA;\cM)=\cM$. We set $C(\cA;\cM) =
\oplus_{n\geq 0}C^n(\cA;\cM)$. Then we introduce the Hochschild
differential $\delta$ on the space $C(\cA;\cM)$ by the formula:
\begin{eqnarray*}
(\delta f)(a_1\otimes\dots \otimes a_{n+1}) &=& a_1 f(a_2\otimes\dots
\otimes a_{n+1}) \nonumber\\ & & + \sum_{i=1}^n (-1)^i f(a_1\otimes
\dots \otimes a_i a_{i+1}\otimes \dots\otimes a_{n+1}) \\ & & +
(-1)^{n+1} f(a_1\otimes \dots \otimes a_n)a_{n+1} \nonumber
\end{eqnarray*}
for any $f\in C^n(\cA;\cM)$. Because $\cA$ is an associative algebra,
one has $\delta^2 = 0$. The cohomology of this differential complex is
denoted by $H(\cA; \cM)$. It is the Hochschild cohomology of $\cA$ with
values in $\cM$.

The bimodule of interest for our purpose is $\cA$ itself. In this case,
the complex $C(\cA;\cA)$ is an associative algebra (See \cite{GER1} and
\cite{GER2} and references therein) and $H(\cA; \cA)$ inherits a
structure of graded commutative algebra.

\bigskip
Let us consider now the previous case with $n=1$. Then $Z^1(\cA;\cA) =
\Im \delta \cap C^1(\cA;\cA)$ is the Lie algebra $\der(\cA)$ of
derivations of $\cA$, and $B^1(\cA;\cA) = \ker \delta \cap C^1(\cA;\cA)$
is the Lie subalgebra $\Int(\cA)$ of $\der(\cA)$ of inner derivations of
$\cA$. This is an ideal of $\der(\cA)$, so $H^1(\cA; \cA)$ is a Lie
algebra, denoted $\out(\cA)$.

\subsection{Relative Hochschild Cohomology}

We follow here the exposition in \cite{GER1} (see also \cite{GER2}). Let
$\cS$ denote a subalgebra of $\cA$, and $\cM$ a bimodule over $\cA$. The
complex $C(\cA,\cS;\cM)$ is defined by
\[ C^0(\cA,\cS;\cM) = \cM^\cS = \{ m\in \cM \,\,/\,\, sm=ms \,\, \forall s\in
\cS \} \]
and $C^n(\cA,\cS;\cM)$ is the linear space of $n$-linear mappings $f :
\cA \otimes \dots \otimes \cA \rightarrow \cM$ such that
\begin{eqnarray*}
f(sa_1\otimes \dots \otimes a_n) &=& s f(a_1\otimes \dots \otimes a_n)
\\ f(a_1\otimes \dots \otimes a_ns) &=& f(a_1\otimes \dots \otimes a_n)s
\\ f(a_1\otimes \dots \otimes a_is\otimes a_{i+1}\otimes \dots\otimes
a_ns) &=& f(a_1\otimes \dots \otimes a_i \otimes sa_{i+1}\otimes
\dots\otimes a_ns)
\end{eqnarray*}
for any $a_i\in\cA$ and $s\in\cS$. $f$ is then a $\cS$-bimodule
homomorphism $\cA \otimes_\cS \dots \otimes_\cS \cA \rightarrow \cM$.

The Hochschild differential $\delta$ maps $C^n(\cA,\cS;\cM)$ into
$C^{n+1}(\cA,\cS;\cM)$, and then defines a cohomology $H(\cA,\cS;\cM)$.
This is the relative Hochschild cohomology of $\cA$ in $\cM$ for $\cS$.
This cohomology can be calculated on a subcomplex of $C(\cA,\cS;\cM)$
\cite{GER1}. Let us denote by $\overline{C}(\cA,\cS;\cM)$ the linear
subspace of $C(\cA,\cS;\cM)$ of elements $f$ such that $f$ vanishes when
at least one of its arguments is in $\cS$. This is the normalized
complex of the relative Hochschild cohomology. These two complexes have
the same cohomology.

\bigskip
Let us now consider the case where $\cS=\cZ(\cA)$. Then the relative
Hochschild cohomology is well adapted to study central bimodules. In
degree 0, one has $C^0(\cA,\cZ(\cA);\cM) = \cM$ for $\cM$ a central
bimodule. In higher degrees, one can remark that $\cA \otimes_{\cZ(\cA)}
\dots \otimes_{\cZ(\cA)} \cA$ is a central bimodule, and then the
normalized relative complex is a set of homomorphisms of central
bimodules.

\bigskip
For future use, consider the case $\cS=\cZ(\cA)$, $\cM = \cA$ and $n=1$.
Then $Z^1(\cA,\cZ(\cA);\cA)$ is exactly the Lie algebra of derivations
of $\cA$ which vanish on the center $\cZ(\cA)$. Remark that
$B^1(\cA,\cZ(\cA);\cA)$ is equal to $B^1(\cA;\cA)$. So, one has the two
left exact sequences:
\[ 0 \rightarrow Z^1(\cA, \cZ(\cA); \cA) \rightarrow \der(\cA) \rightarrow
\der(\cZ(\cA)) \]
\[ 0 \rightarrow H^1(\cA, \cZ(\cA); \cA) \rightarrow \out(\cA) \rightarrow
\der(\cZ(\cA)) \]
which are not short exact sequences in general. The condition
$H^1(\cA,\cZ(\cA);\cA)=0$, which means that any derivation of $\cA$
which vanishes on $\cZ(\cA)$ is an inner derivation, gives the
injectivity of the canonical homomorphism $\out(\cA) \rightarrow
\der(\cZ(\cA))$.

\subsection{Constrained Hochschild Cohomology}

Let us now introduce a new subcomplex of the Hochschild complex. As
before, $\cA$ is an associative algebra with unit and $\cM$ is a
bimodule over $\cA$. Let $\cC$ be an ideal in $\cA$ and $\cN$ a
sub-bimodule of $\cM$ such that $cm, mc\in\cN$ for any $c\in\cC$ and
$m\in\cM$. This is equivalent to say that $\cC$ is included in the two
side ideal
\[ \cI_\cN = \{ a\in\cA \,\, /\,\, a\cM \subset \cN \mbox{ and } \cM a \subset
\cN \} \]

\bigskip
We define the subcomplex $C(\cA,\cC;\cM,\cN)$ of $C(\cA;\cM)$ of the
mappings $f : \cA\otimes \dots \otimes\cA \rightarrow \cM$ such that
$f(a_1\otimes\dots\otimes a_n)\in\cN$ if at least one of the $a_i$ is in
$\cC$. In degree 0, $C^0(\cA,\cC;\cM,\cN) = \cM$. It is easy to see that
this subcomplex is stable by the Hochschild differential $\delta$. So
one has a cohomology $H(\cA,\cC;\cM,\cN)$. This is the constrained
cohomology of $\cA$ in $\cM$ by $(\cC,\cN)$. One has then the following
Lemma:
\begin{lemme}
\label{cohom}
In the above situation, one has a canonical mapping of graded vector
spaces
\[ H(\cA,\cC;\cM,\cN) \rightarrow H(\cA/\cC;\cM/\cN) \]
where the second cohomology is the ordinary Hochschild cohomology of the
bimodule $\cM/\cN$ over the algebra $\cA/\cC$.
\end{lemme}

\medskip
{\sl Proof:} Let $pr : \cM \rightarrow \cM/\cN$ denote the projection
from the bimodule $\cM$ over $\cA$ on the bimodule $\cM/\cN$ over
$\cA/\cC$, and $a\rightarrow [a]$ the projection $\cA\rightarrow
\cA/\cC$. Then one has $pr(am) = [a]pr(m)$ for any $a\in\cA$ and $m\in
\cM$, and a similar formula for $ma$.

Any $f\in C(\cA,\cC;\cM,\cN)$ can be mapped into $\chi(f)\in
C(\cA/\cC;\cM/\cN)$ by the definition
\[ \chi(f)([a_1]\otimes \dots \otimes [a_n]) =
(pr\circ f)(a_1 \otimes \dots \otimes a_n) \]
Then it is easy to see that
\[ pr\circ \delta = \overline{\delta} \circ pr \]
where $\delta$ is the Hochschild differential on $C(\cA,\cC;\cM,\cN)$
and $\overline{\delta}$ the Hochschild differential on
$C(\cA/\cC;\cM/\cN)$.~\hfill$\square$

\bigskip
A simpler situation occurs when one takes $\cM=\cA$ and $\cN=\cC$. Then
the subcomplex is a subalgebra of $C(\cA;\cA)$, but not an ideal. We
denote it by $C_\cC(\cA;\cA)$, and its cohomology by $H_\cC(\cA;\cA)$.
In degree 1, one has obviously $B^1_\cC(\cA;\cA) = B^1(\cA;\cA)$.
$Z^1_\cC(\cA;\cA)$ is the Lie algebra of derivations of $\cA$ which
preserve $\cC$. $B^1_\cC(\cA;\cA)$ is an ideal in this Lie algebra. Then
$H^1_\cC(\cA;\cA)$ is a Lie algebra.

\section{Noncommutative Submanifolds}

In this section, we introduce a noncommutative generalization of the
notion of submanifold of a manifold.

\subsection{The Commutative Case}

We first recall the situation in the commutative case. Let $M$ be a
smooth compact manifold, and let $N\subset M$ be a closed submanifold.
Any smooth function $f:M \rightarrow \gR$ can be restricted to $N$. Thus
one has a mapping
\[ \cF(M) \buildrel{p}\over\rightarrow \cF(N) \]
where $\cF(M)$ is the algebra of smooth functions on $M$.  This mapping
is in fact surjective, and there exists a short exact sequence
\[ 0\rightarrow \cC \rightarrow \cF(M) \buildrel{p}\over\rightarrow
 \cF(N) \rightarrow 0 \]
where $\cC$ is the ideal of $\cF(M)$ of functions vanishing on $N$.

\bigskip
A vector field $X\in\Gamma(M)$ on $M$, which satisfies $Xf\in \cC$ for
any $f\in \cC$, can be restricted to a vector field $\overline{X}$ on
$N$. Thus one has an homomorphism of Lie algebras
\[ \Gamma_\cC(M) \buildrel{\pi}\over\rightarrow \Gamma(N) \]
where $\Gamma_\cC(M) = \{ X\in\Gamma(M)\; /\; X\cC\subset\cC \}$. This
mapping is surjective, and there exists a short exact sequence of Lie
algebras:
\[ 0 \rightarrow \Gamma_\cF \rightarrow \Gamma_\cC(M)\buildrel{\pi}\over
\rightarrow \Gamma(N) \rightarrow 0 \]
where $\Gamma_\cF = \{ X\in\Gamma(M)\; /\; X\cF(M)\subset\cC \}$ is an
ideal of the Lie algebra $\Gamma_\cC(M)$.

\subsection{The Noncommutative Case}

Now we can generalize these notions to the framework of noncommutative
geometry. Let $\cA$ be an associative algebra over $\gC$ with unit and
let $\cC$ be an ideal in $\cA$. We denote by $\cQ = \cA/\cC$ the
quotient algebra and $p: \cA \rightarrow \cQ$ the quotient mapping.

We can consider the two following Lie subalgebras of $\der(\cA)$:
\[ \cG_\cC = \{ X\in \der(\cA)\; /\; X\cC\subset \cC \} \]
and
\[ \cG_\cA = \{ X\in \der(\cA) \; /\; X\cA\subset \cC \} \]
One sees that $\cG_\cA$ is an ideal in $\cG_\cC$. One has a mapping
$\cG_\cC \buildrel{\pi}\over\rightarrow \der(\cQ)$ defined by $\pi(X)
p(a) = p(Xa)$ for any $a\in\cA$ and $X\in\cG_\cC$. The kernel of this
mapping is exactly $\cG_\cA$.

\bigskip
\begin{definition}
\label{defsub}
The quotient algebra $\cQ = \cA/\cC$ will be called a {\sl submanifold
algebra of $\cA$} if $\pi$ is surjective. The ideal $\cC$ of $\cA$ is
called the {\sl constraint ideal} for $\cQ$.
\end{definition}

In this situation, one has the short exact sequence of Lie algebras
\begin{equation}
\label{secdersub}
 0 \rightarrow \cG_\cA \rightarrow \cG_\cC
\buildrel{\pi}\over\rightarrow \der(\cQ) \rightarrow 0 \end{equation}
The condition of the definition imposes a strong relation between the
differential structure on $\cA$ and the differential structure on $\cQ$.
This strong relation is revealed in the following Proposition:

\begin{proposition}
There exists a short exact sequence of graded differential algebras
\begin{equation}\label{secformes}
0 \rightarrow \Omega_{\der,\cC} \rightarrow \Omega_\der(\cA)
\buildrel{p}\over\rightarrow \Omega_\der(\cQ) \rightarrow 0
\end{equation}
\end{proposition}

\medskip
{\sl Proof:} Let $\overline{X}=\pi(X)\in \der(\cQ)$ for any
$X\in\cG_\cC$ and let $\overline{a} = p(a)\in \cQ$ for any $a\in \cA$.
Then one has over $\cQ$
\[ d \overline{a}(\overline{X}) = \overline{X}\overline{a} = p(Xa) = p(da(X))
\]
One then extends $p$ in a mapping $\Omega_\der^n(\cA) \rightarrow
\Omega_\der^n(\cQ)$ by the relation
\[ p(a_0da_1 \dots da_n) = p(a_0) dp(a_1) \dots dp(a_n) \]
and then one has
\[ d\circ p = p\circ d \]
and
\[ i_{\overline{X}} \circ p = p\circ i_X \]
It is easy to see that $p$ is surjective; so we obtain the short exact
sequence (\ref{secformes}).~\hfill$\square$

\bigskip
{\bf Remarks:}

1. In the short exact sequence (\ref{secformes}), one has
\[ \Omega_{\der,\cC}^n = \{ \omega\in \Omega_\der^n(\cA) \;/\;
\forall X\in\cG_\cC, \; i_X \omega\in \Omega_{\der,\cC}^{n-1} \} \]
with $\Omega_{\der,\cC}^0 = \cC$. For example, for any $a\in\cC$,
$da\in\Omega_{\der,\cC}^1$.

\medskip
2. Nothing can be said about any canonical relation between
$\underline{\Omega}_\der(\cA)$ and $\underline{\Omega}_\der(\cQ)$.

\bigskip
Let us now study the derivations of $\cQ$. Any inner derivation of $\cA$
is obviously in $\cG_\cC$. In the quotient homomorphism $\cG_\cC
\buildrel{\pi}\over\rightarrow \der(\cQ)$, these inner derivations are
mapped on inner derivations, from the very definition of $\pi$. It is
easy to see that $\pi$ restricted to inner derivations is surjective on
inner derivations of $\cQ$ (even if $\pi$ does not satisfy the condition
of Definition~\ref{defsub}, {\sl i.e.} $\pi$ is not surjective) and one
has $\pi(ad(a)) = ad( p(a))$ for any $a\in \cA$. So, the kernel of $\pi$
contains $ad(\cC) =\{ ad(c) / c\in\cC\} \subset \der(\cA)$.

\begin{lemme}
\label{lemmederint}
If $\cQ = \cA / \cC$ has only inner derivations, then the mapping
$\cG_\cC \rightarrow \der(\cQ)$ is surjective. Then $\cQ$ is a
submanifold algebra.
\end{lemme}

\medskip
{\sl Proof:} This is a direct consequence of the previous discussion
about inner derivations.~\hfill$\square$

\bigskip
It is now interesting to say something about the other derivations of
$\cQ$, that is, about the first Hochschild cohomology of $\cQ$ with
values in itself. One has the Lemma:

\begin{lemme}
One has a surjective homomorphism of Lie algebras
\[ H_\cC^1(\cA;\cA) \rightarrow H^1(\cQ;\cQ) \]
\end{lemme}

\medskip
{\sl Proof:} This is a direct consequence of Lemma~\ref{cohom} (with
$\cM=\cA$ and $\cN = \cC$), the previous remark about inner derivations,
and the surjectivity of $\pi$ from
Definition~\ref{defsub}.~\hfill$\square$

\bigskip
One can say something about the kernel of this mapping, if one imposes a
supplementary condition on the ideal $\cC$.

\begin{proposition}
If the constraint ideal $\cC$ for the submanifold algebra $\cQ$
satisfies
\[ ad(\cC) = \{ ad(a) \,\,/\,\, a\in\cA \mbox{ and } [a,\cA]\subset\cC \} \]
or equivalently, if $\ker\pi\cap \Int(\cA) = ad(\cC)$, then one has the
short exact sequence of Lie algebras
\begin{equation}
\label{seccohom}
 0 \rightarrow H^1(\cA;\cC) \rightarrow H_\cC^1(\cA;\cA) \rightarrow
H^1(\cQ;\cQ)\rightarrow 0
\end{equation}
\end{proposition}
In $H^1(\cA;\cC)$, $\cC$ is considered as a bimodule over $\cA$.

\medskip
{\sl Proof:} The condition on $\cC$ means in fact that one has the short
exact sequence of Lie algebras
\[ 0 \rightarrow B^1(\cA;\cC) \rightarrow B_\cC^1(\cA;\cA)=B^1(\cA;\cA)
\rightarrow B^1(\cQ;\cQ)\rightarrow 0 \]
The new information is the exactness at $B_\cC^1(\cA;\cA)$. If one
associates this short exact sequence with the short exact
sequence~(\ref{secdersub}) written as
\[ 0 \rightarrow Z^1(\cA;\cC) \rightarrow Z_\cC^1(\cA;\cA) \rightarrow
Z^1(\cQ;\cQ)\rightarrow 0 \]
then one obtains the exactness of (\ref{seccohom}).~\hfill$\square$

\bigskip
In algebraic geometry (\cite{SCH1} and references therein), one works
with the commutative algebra with unit $\cA=\gC[X_1, \dots, X_n]$ of
complex polynomials of $n$ variables. The geometric objects are
considered as zero sets of polynomials. An ideal $\cC$ represents the
set of points $V(\cC) = \{ x\in\gC^n \,\, /\,\, P(x) = 0 \,\, \forall
P\in\cC \}$. From the point of view of the ``duality'' {\sl set of
points} $\leftrightarrow$ {\sl algebra of functions}, the set $V(\cC)$
is represented by the algebra $\cQ=\gC[X_1, \dots, X_n]/\cC$. If $\cQ$
admits ideals, then the set $V(\cC)$ admits subsets. But if $\cQ$ does
not have any ideal, then the set $V(\cC)$ can be considered as a {\sl
point}. This is equivalent to the fact that $\cC$ is a maximal ideal in
$\cA$. Notice that from the point of view of ordinary geometry, points
are the minimal sets. The only maximal ideals of $\cA$ are generated by
$n$ polynomials $X_i - a_i$ where $a_i\in\gC$. The point represented by
this ideal is obviously $(a_i)\in\gC^n$. Notice that maximal ideals are
in one-to-one correspondence with the caracters of $\cA$. The quotient
mapping $\cA\rightarrow \cA/\cC$ is the restriction at the set of points
represented by $\cC$. If $\cC$ is maximal, the restriction of an element
$P\in\cA$ at $\cC$ is exactly the value of this polynomial at the point
of $\gC^n$ represented by $\cC$.

This correspondence {\sl points} $\leftrightarrow$ {\sl maximal ideals}
is also used in the theory of commutative Banach algebras and
commutative $C^\ast$-algebras \cite{LAN1}. In this context, maximal
ideals are also in one-to-one correspondence with characters. If $\cC$
is a maximal ideal in the commutative Banach algebra with unit $\cA$,
then the quotient $\cA/\cC$ is isomorphic to $\gC$. (The quotient
mapping $\cA\rightarrow \cA/\cC$ is a character). In the theory of
commutative $C^\ast$-algebras, by the Gel'fand transformation, the set
of characters is exactly the set of points, on which it is possible to
put a canonical topology. One says that a point takes its values in the
quotient $\cA/\cC\simeq\gC$. So, in those two situations, points are
maximal ideals, and take their values in the field $\gC$.

In noncommutative geometry, an ideal $\cC$ of a given complex algebra
with unit $\cA$ can also be interpreted as a ``subspace'' of the non
commutative ``space'' dualy represented by $\cA$. This subspace can be
considered as a ``submanifold'' if the differential structure of
$\cA/\cC$ is compatible with the differential structure of $\cA$. One of
the compatibility conditions one can take is Definition~\ref{defsub}.

Now, if the ideal is maximal, then the quotient algebra is simple. It is
then a ``point'', in the sense that it can not have ``subspace''. But
then, considering the quotient $\cQ=\cA/\cC$, one sees that points take
their values in ({\sl a priori} noncommutative) simple algebras, and not
in fields as in the commutative case. There is then a residual
structure, of purely noncommutative origin. See Example~5 below for
applications in physics.

To any ideal $\cC$ in $\cA$, one can construct
$\cG_\cA\subset\der(\cA)$. If $\cC$ is a maximal ideal, the quotient of
linear space $T_\cC = \der(\cA)/\cG_\cA$ can be considered as the
``tangent space'' at the point $\cC$ in the ``manifold'' represented by
$\cA$. The value of a derivation $X$ at the ``point'' $\cC$ is the image
of $X$ by the quotient mapping $\der(\cA)\rightarrow T_\cC$. One can
also take the value of a 1-form $\alpha\in\Omega^1(\cA)$ at $\cC$ by the
definition $\alpha_\cC :T_\cC \rightarrow \cQ$, $\alpha_\cC(X_\cC) =
p\circ \alpha(X)$ for any $X\in\der(\cA)$ whose value at $\cC$ is
$X_\cC$. This can be generalized for any $n$-form in $\Omega^n(\cA)$.

\subsection{Examples}

\hspace*{\parindent}{\bf Example 1:} The commutative case.

In the commutative case, any smooth closed submanifold of a smooth
compact manifold gives a submanifold algebra: the algebra of smooth
functions on this submanifold.

\bigskip
{\bf Example 2:} The tensor algebra.

Let $\cA$ be the free algebra with unit over $\gC$ generated by $n$
elements $x^1, \dots, x^n$, with $n\geq 2$.

Any derivation of $\cA$ is given by $n$ elements $P^i(x^1, \dots ,
x^n)\in\cA$. We denote it by $D = (P^i)_{1\leq i\leq n}$. The value of
this derivation on any element of $\cA$ is obtained by the Leibniz rule
and the definition $D(x^i) = P^i(x^1, \dots , x^n)$.

If one takes $\cC$ the ideal in $\cA$ generated by $x^1$ then the
algebra $\cQ$ is the free algebra with unit over $\gC$ generated by
$x^2, \dots , x^n$, and one has $\cA = \cC\oplus\cQ$.

Any derivation in $\cG_\cC$ is the sum of two kinds of derivation:
$(P^i)_{1\leq i\leq n}$, with $P^i\in\cQ$ and $P^1=0$ ($\cQ$ is
considered here as a subalgebra of $\cA$), and $(P^i)_{1\leq i\leq n}$,
with $P^i\in\cC$. Any derivation in $\cG_\cA$ is of the second kind, and
the Lie algebra of derivation of $\cQ$ is the set of the first kind
derivations in $\cG_\cC$. So the condition of the
Definition~\ref{defsub} is fulfilled. $\cQ$ is thus a submanifold
algebra of $\cA$. In this case, one has $\cG_\cC =
\cG_\cA\oplus\der(\cQ)$.

As maximal ideals of $\cA$, one has the ideals generated by the $n$
elements $x^i - a_i$ where $a_i\in\gC$. Then the point associated to
such an ideal is a point in $\gC^n$, with values in $\gC$. This
situation is analogous to the situation of the polynomial algebra
generated by the $n$ variables $x^i$, for which there are only those
maximal ideals. It is not difficult to see that such an ideal contains
the ideals generated by the expressions $x^ix^j - x^jx^i$. In the case
of the tensor algebra $\cA$, there are other interesting maximal ideals,
as the following examples show.

\bigskip
{\bf Example 3:} The Heisenberg algebra.

Let $\cA$ be the free algebra with unit generated by two elements $x,y$.
Consider in $\cA$ the ideal generated by $xy-yx-i\bbbone$. Then the
quotient algebra is the Heisenberg algebra $\cH$, generated by two
elements $p,q$ and the relation $pq-qp=i\bbbone$. It is well known that
this algebra is simple. The ideal is maximal. In the quotient, we take
$x\mapsto p$ and $y\mapsto q$.

Now let us consider derivations. If we denote by $D= (X,Y)$ the
derivation $D(x) = X$ and $D(y)=Y$, then one has
\[ \cG_\cC = \left\{ (X,Y) \,\, / \,\, [X,y]+[x,Y]\in\cC \right\} \]
where $X,Y\in\cA$, and
\[ \cG_\cA = \left\{ (X,Y) \,\, / \,\, X,Y\in\cC \right\} \]
On the other hand, one knows that $\cH$ has only inner derivations (See
\cite{MDV2}, for instance.), so
\[ \der(\cH) = \left\{ ([A,p],[A,q]) \,\, / \,\, A\in\cH \right\} \]
with the same notations as above. It is easy to prove that the mapping
$\cG_\cC \rightarrow \der(\cH)$ (the quotient by $\cG_\cA$) is
surjective (one can use Lemma~\ref{lemmederint}, but the direct
calculation shows in this particular case how that works). Indeed, take
$A\in\cH$ and let $\overline{A}\in\cA$ be such that $\overline{A}
\mapsto A$ in the quotient mapping $\cA \rightarrow \cH$. Then the
derivation $([\overline{A} , x], [\overline{A}, y])$ maps to
$([A,p],[A,q])\in\der(\cH)$. One must then show that this derivation is
indeed in $\cG_\cC$. This is equivalent to showing $[[\overline{A} , x],
y] + [x, [\overline{A}, y]] \in\cC$. But this expression equals
$-[[x,y], \overline{A}]$ which is obviously in the kernel of the mapping
$\cA \rightarrow \cH$. So it is in $\cC$.

The Heisenberg algebra is then a submanifold algebra, which can be
regarded, from the point of view of algebraic geometry, as a point in
the free algebra with unit $\cA$. Its tangent space is the linear space
$\cH\oplus\cH$.

\bigskip
{\bf Example 4:} The matrix algebra.

Let $\cA$ denote as above the free algebra with unit generated by two
elements $x,y$. Let $q\in \gC$ a $n^{\mathrm th}$ unit root, $q^n = 1$.
Let $\cC$ denote the ideal in $\cA$ generated by the relations
\[ xy - qyx, \quad x^n - \bbbone, \quad y^n - \bbbone \]
and denote by $U$ and $V$ the images of $x$ and $y$ in the quotient
mapping $\cA \rightarrow \cQ$. Let us show that this algebra is the
matrix algebra $M(n,\gC)$. Any element of $\cQ$ can be written as
\[ \sum_{0\le k,\ell \le n-1} a_{k,\ell} U^k V^\ell \]
so $\dim\cQ\le n^2$. Now, the following two matrices in $M(n,\gC)$,
\[
U = \left(
\begin{array}{ccccc}
1 & 0 & 0 & \cdots & 0 \\ 0 & q & 0 & \cdots & 0 \\ 0 & 0 & q^2 & \cdots
& 0 \\
\vdots & \vdots &  & \ddots & \vdots \\
0 & 0 & 0 & \cdots & q^{n-1}
\end{array}
\right), \quad
V = \left(
\begin{array}{ccccc}
0 & 1 & 0 & \cdots & 0 \\ 0 & 0 & 1 & \cdots & 0 \\
\vdots & \vdots & \ddots & \ddots & \vdots \\
0 & 0 & \cdots & 0 & 1 \\ 1 & 0 & \cdots & 0 & 0
\end{array}
\right)
\]
satisfy the relations of the algebra $\cQ$, and then generate a
subalgebra of $M(n,\gC)$. Because the only matrices which commute with
this subalgebra are the multiple of identity, this is the full matrix
algebra.

It is well known that the matrix algebra has only inner derivations. By
Lemma~\ref{lemmederint}, $M(n,\gC)$ can be considered as a submanifold
algebra of the tensor algebra. Notice that this algebra is simple, and
so can be considered as a ``point'' in the tensor algebra.

\bigskip
{\bf Example 5:} The matrix value functions.

Consider, as in \cite{DKM2}, the algebra $\cA = C^\infty(V)\otimes
M(n,\gC)$ of matrix value functions on a manifold $V$. Let $p\in V$ any
point of the manifold. Take $\cC$ the ideal of functions vanishing at
$p$. This is obviously a maximal ideal. It has been shown in \cite{DKM2}
that $\der(\cA) = [\der(C^\infty(V))\otimes\bbbone] \oplus [ C^\infty(V)
\otimes \der(M(n,\gC))]$. Then a simple calculation shows that $\cQ =
\cA / \cC$ is the matrix algebra $M(n,\gC)$, and is a submanifold
algebra of $\cA$. The ``tangent space'' at $p$ is $T_pV \oplus
\der(M(n,\gC))$, where $T_pV$ is the ordinary tangent space of $V$ at
$p$.

The physical interpretation of this situation is the following: from the
point of view of noncommutative differential geometry, each point of
space-time is a matrix, instead of $\gC$ (or $\gR$) in ordinary
differential geometry. The structure looks like a fiber bundle, the
fiber been a matrix algebra, but the differential structure is
different, because the purely noncommutative differential structure of
the matrix algebra (which is far from being trivial) is taken into
account at each point of $V$. This supplementary differential structure
of points has important consequences for gauge fields theory, as has
been shown in \cite{DKM2}.

This situation can be modified without many changes, by taking the
algebra of sections of bundle over $M$, with fiber $M(n,\gC)$.

\section{Noncommutative Quotient Manifolds}

In this section, we introduce a generalization to the noncommutative
framework of the notion of quotient manifold. We then introduce the
generalization of the action of a group on a manifold which gives a way
to construct such quotient manifolds. We give examples and we examine
the possible relations with connections on central bimodules.

\subsection{Quotient Manifold Algebra}

Let $\cA$ be an associative algebra with unit. Let $\cB$ be a subalgebra
of $\cA$. Then we define the Lie subalgebras of $\der(\cA)$:
\[ \hat{\fg} = \{ X\in\der(\cA) \; / \; X\cB=0 \} \]
and
\[ \fh = \{ X\in\der(\cA) \; / \; X\cB\subset \cB \} \]
Notice that $\hat{\fg}$ is an ideal in $\fh$, {\sl i.e.} $[\fh,
\hat{\fg}] \subset \hat{\fg}$.

\bigskip
One has a natural homomorphism of Lie algebras $\rho : \fh \rightarrow
\der(\cB)$, $X\mapsto \tilde{X}$, the restriction of $X$ to $\cB$. The
kernel of this homomorphism is exactly $\hat{\fg}$.

\begin{definition}
\label{defquo}
The subalgebra $\cB$ of $\cA$ is a {\sl quotient manifold algebra of
$\cA$} if the following three conditions are fulfilled

$(i)$ $\cZ(\cB) = \cB \cap \cZ(\cA)$,

$(ii)$ $\der(\cB) \simeq \fh / \hat{\fg}$,

$(iii)$ $\cB = \{ a\in\cA \,\, / \,\, Xa=0 \,\, \forall X\in\hat{\fg}
\}$.
\end{definition}

\bigskip
Condition $(i)$ gives to $\fh$ and $\hat{\fg}$ a structure of
$\cZ(\cB)$-module. $\fh / \hat{\fg}$ is then naturally a
$\cZ(\cB)$-module, and condition $(ii)$ is an isomorphism of
$\cZ(\cB)$-modules. One then has the short exact sequence of Lie
algebras and $\cZ(\cB)$-modules
\begin{equation}
\label{secderquo}
0 \rightarrow \hat{\fg} \rightarrow \fh \buildrel{\rho}\over\rightarrow
\der(\cB) \rightarrow 0
\end{equation}
Now, the Lie subalgebra $\hat{\fg}$ of $\der(\cA)$ gives a Cartan
operation on $\underline{\Omega}_\der(\cA)$. Condition $(iii)$ says that
$\cB$ is exactly the basic algebra in $\cA$ for this operation.

Let $\omega \in \underline{\Omega}_\der^n(\cA)$ be a basic element for
the operation of $\hat{\fg}$, $i_X\omega = 0$, and $L_X \omega=0$ for
any $X\in \hat{\fg}$. Then $d\omega$ is also basic. One can then define
$\tilde{\omega} \in \underline{\Omega}_\der^n(\cB)$ by the relation
\[ \tilde{\omega}(\tilde{X}_1, \dots, \tilde{X}_n) = \omega(X_1, \dots , X_n)
\]
for any $\tilde{X}_1, \dots, \tilde{X}_n \in\der(\cB)$ and any
representatives $X_1, \dots , X_n\in \fh$. By the Koszul formula and
condition $(iii)$, it is easy to show that $\omega(X_1, \dots ,
X_n)\in\cB$ for $X_i\in\fh$. Note that condition $(ii)$ is essential to
ensure the consistency of this definition. Condition $(i)$ shows that
$\tilde{\omega}$ is $\cZ(\cB)$-linear. The Koszul formula shows then
that $d\tilde{\omega} = \widetilde{d\omega}$.

So one has the Lemma:
\begin{lemme}
One has a mapping of graded differential algebras
\[ \underline{\Omega}_{\der,B}(\cA) \rightarrow \underline{\Omega}_\der(\cB) \]
where $\underline{\Omega}_{\der,B}(\cA)$ is the sub-algebra of
$\underline{\Omega}_\der(\cA)$ of basic elements for $\hat{\fg}$.
\end{lemme}

\bigskip
{\bf Remarks:}

1. In degree 0, by the very definition, one has
$\underline{\Omega}_{\der,B}^0(\cA) = \cB =
\underline{\Omega}_\der^0(\cB)$.

\medskip
2. No canonical mapping can be constructed between the basic elements of
$\Omega_\der(\cA)$ and $\Omega_\der(\cB)$ without more information on
the algebras $\cA$ and $\cB$.

\medskip
3. Condition $(ii)$ can be relaxed if we define $\der(\cB)$ to be the
Lie algebra $\fh / \hat{\fg}$, even if $\cB$ accepted other derivations.
In this situation, one has a kind of induced differential structure on
$\cB$ (see example 1 below).

\begin{proposition}
\label{isom}
If the $\cZ(\cA)$-module induced by $\fh$ in $\der(\cA)$ is $\der(\cA)$
itself, then we have an isomorphism of graded differential algebras
\[ \underline{\Omega}_{\der, B}(\cA) \simeq \underline{\Omega}_\der(\cB). \]
\end{proposition}

\medskip
{\sl Proof:} First, let us prove that the mapping
$\underline{\Omega}_{\der,B}(\cA) \rightarrow
\underline{\Omega}_\der(\cB)$ constructed above is injective. If
$\tilde{\omega}$ is zero in $\underline{\Omega}_\der^n(\cB)$ then for
any $X_1, \dots , X_n\in \fh$ we have $\omega(X_1, \dots , X_n)=0$. Now,
$\omega$ is $\cZ(\cA)$-linear, so $\omega$ is zero on the
$\cZ(\cA)$-module induced by $\fh$ in $\der(\cA)$. This proves
injectivity.

Let $\tilde{\omega}\in \underline{\Omega}_\der^n(\cB)$ be any $n$-form.
Define $\omega$ as an antisymmetric $n$-$\cZ(\cB)$-linear mapping from
$\fh\otimes_{\cZ(\cB)}\dots \otimes_{\cZ(\cB)}\fh$ to $\cB$ by the
relation
\[ \omega(X_1, \dots , X_n) = \tilde{\omega}(\tilde{X}_1, \dots,
\tilde{X}_n) \in\cB\subset\cA \]
for any $X_1, \dots , X_n\in \fh$. Then we extend $\omega$ on the
$\cZ(\cA)$-module induced by $\fh$, by $\cZ(\cA)$ linearity. Notice that
$\omega$ is already $\cZ(\cB)$ linear. By hypothesis, $\omega$ is then
an element of $\underline{\Omega}_\der^n(\cA)$. We have $i_X \omega = 0$
for any $X\in \hat{\fg}$, so $\omega$ is horizontal for the action of
$\hat{\fg}$ in $\underline{\Omega}_\der(\cA)$. Now, notice that the
$(n+1)$-form in $\underline{\Omega}_\der(\cA)$ which comes from
$d\tilde{\omega}$ is exactly $d\omega$, because by the Koszul formula
they coincide on $\fh$. So $d\omega$ is also horizontal, and then
$\omega$ is basic in $\underline{\Omega}_\der(\cA)$. This proves
surjectivity.~\hfill$\square$

\subsection{Action}

Let $M$ be a manifold and $G$ a Lie group. An action of $G$ on $M$ gives
a Lie algebra homomorphism $\fg \rightarrow \Gamma(M)$ from the Lie
algebra $\fg$ of $G$ to the Lie algebra of vector fields on $M$. Then
one has an Cartan operation of the Lie algebra $\fg$ on the graded
commutative differential algebra $\Omega(M)$ of de~Rham differential
forms on $M$.

In the noncommutative case, we will say we have an action of the Lie
algebra $\fg$ on an associative algebra with unit $\cA$ if there is a
homomorphism of Lie algebras $\fg\rightarrow \der(\cA)$. In this
situation one has an operation of $\fg$ on the graded differential
algebra $\underline{\Omega}_\der(\cA)$.

Then one can take as subalgebra of $\cA$ the basic algebra for this
operation. In this situation, if $\cB$ is a quotient manifold algebra,
then one has the noncommutative version of the quotient manifold by the
``leaves'' defined by $\fg$.

If the homomorphism $\fg\rightarrow \der(\cA)$ is injective (take the
image of $\fg$ if not), one can identify $\fg$ with its image. Then, one
has the inclusion $\fg\subset \hat{\fg}$, but the equality is not the
generic case. Between these two Lie algebras, one has a third one, the
$\cZ(\cA)$-module induced by $\fg$ in $\der(\cA)$, denoted by
$\fg_{\cZ(\cA)}$. If $\cB$ is the basic algebra in $\cA$ for the
operation of $\fg$, the condition $(iii)$ of Definition~\ref{defquo} is
fulfilled.

\subsection{Examples}

\hspace*{\parindent}{\bf Example 1:} The inner derivations.

Let $\cA$ be an associative algebra with unit for which there are inner
derivations. Suppose one has $H^1(\cA,\cZ(\cA);\cA) = 0$. Take the
operation of $\fg=\Int(\cA)$ on $\cA$. Then one has $\cB= \cZ(\cA)$,
$\hat{\fg} = \fg$ because the first relative cohomology group vanishes,
and $\fh = \der(\cA)$. Take then the induced differential structure on
$\cB$ by setting $\der(\cB) = \out(\cA)= \fh/\fg$. The algebra of
differential forms associated to $\cB$ is then, by
Proposition~\ref{isom}, the algebra $\underline{\Omega}_\out(\cA)$
introduced in \cite{DVM2}.

\bigskip
{\bf Example 2:} The noncommutative torus.

\nopagebreak
Let $T_q$ denote the complex associative algebra with unit of elements
of the form
\[ a= \sum_{k,\ell\in\gZ} c_{k\ell}U^k V^\ell \]
with
\[ \| a \|_m = \sup_{k,\ell\in\gZ} |c_{k\ell}|  (1+|k| + |\ell|)^m <\infty \]
and the relation
\[ UV = qVU \]
for $q\in\gC$ such that $q^N=1$ for $N\in\gN$. We take $N$ the minimal
one for which this is true. The center of this algebra is the set of
elements depending only of $U^N$ and $V^N$.

The derivations of this algebra are the inner derivations and the
derivations of the form
\[ a(U^N,V^N)D_U + b(U^N,V^N)D_V \]
where $D_U(U) = U$, $D_U(V)=0$, $D_V(U)=0$ and $D_V(V)=V$, and
$a(U^N,V^N)$ and $b(U^N,V^N)$ belong to $\cZ(T_q)$.

Take $\fg=\Int(T_q)$, then $\cB=\cZ(T_q)$, $\hat{\fg} = \fg$ and $\fh =
\der(T_q)$. We are then in the situation of the previous example. Then
the differential algebra of forms on $\cZ(T_q)$ is the basic algebra of
the differential algebra of forms on $T_q$. But now, remark that the
center $\cZ(T_q)$ is isomorphic to the algebra $C^\infty(S^1\times S^1)$
of smooth functions on the (ordinary) torus. This isomorphism is
$U^N\mapsto e^{2\pi it}$ and $V^N\mapsto e^{2\pi is}$. Then an element
$a\in\cZ(T_q)$ is mapped on the Fourier expansion of an element of
$C^\infty(S^1\times S^1)$. Thus, the algebra of forms on $\cZ(T_q)$ is
the de~Rham algebra of forms on the torus.

\subsection{Connections}

Let $\cB$ be a quotient manifold algebra of $\cA$. Then $\cA$ is a
central bimodule over the algebra $\cB$. Let $\psi :
\der(\cB)\rightarrow \fh$ be a spliting
of the
short exact sequence (\ref{secderquo}), considered as a short exact
sequence of $\cZ(\cB)$-modules (forgetting the Lie algebra structures).

\begin{proposition}
\label{connection}
For any $X\in\der(\cB)$, the mapping
\begin{eqnarray*}
\nabla_X : \cA &\rightarrow & \cA \\
a &\mapsto & \psi(X)a
\end{eqnarray*}
is a connection on the central bimodule $\cA$ over $\cB$.

The curvature of this connection is the obstruction on $\psi$ to be a
spliting of the short exact sequence (\ref{secderquo}) of Lie algebras.
\end{proposition}

\medskip
{\sl Proof:} This is an immediate consequence of the fact that $\psi$ is
a $\cZ(\cB)$-modules homomorphism, such that $\psi(X)b = Xb$ for any
$b\in\cB\subset\cA$. The curvature of this connection is
\[ R(X,Y) = [\psi(X), \psi(Y)] - \psi([X,Y]) \]
which proves the Proposition.~\hfill$\square$

\bigskip
Such a connection gives a projection $P: \fh \rightarrow
\hat{\fg}\subset \fh$ of $\cZ(\cB)$-modules defined by $P(X) = X
-\psi\circ\rho(X)$. Then one has $\fh = \ker P \oplus \hat{\fg}$.

Conversely, a projection $P: \fh \rightarrow \hat{\fg}\subset\fh$ of
$\cZ(\cB)$-modules defines a split, and so a connection on $\cA$.

\bigskip
Let $P(M,G)$ be a principal bundle, where $M$ is the base manifold and
$G$ the structure group, and let $\fg$ be its Lie algebra. Denote by
$\cA$ the (commutative) algebra of smooth functions on $P$, and $\cB$
the algebra of smooth functions on $M$. Then one can consider that
$\cB\subset\cA$ because of the projection $P\rightarrow M$. The Lie
algebra $\fg$ can be injectively mapped into $\Gamma(P)$, the vector
fields on $P$, and more precisely, into the vertical vector fields. Thus
$\fg$ operates on $\cA$. The algebra $\cB$ is obviously the basic
algebra for this operation and $\hat{\fg}$ is exactly the Lie algebra of
vertical vector fields on $P$.

It is well known that a connection on $P$ can be given as a $\cB$-linear
mapping $\Gamma(M) \rightarrow \Gamma(P)$, $X\mapsto X^h$, the
horizontal lift, with its usual properties, one of them being $[\fg,
X^h] = 0$ for all $X\in\Gamma(M)$. In fact this mapping is a spliting of
(\ref{secderquo}) (remember here that $\cB=\cZ(\cB)$).

Then one could think that the connections introduced by the construction
of Proposition~\ref{connection} are generalizations of connections on
principal bundles. But this is not completely true, because a principal
bundle has many more properties than a couple $(\cA, \cB)$ of an algebra
and a quotient manifold algebra. For example, one can introduce
associated bundles, on which connections can be transported.

In order to obtain a similar situation, one must introduce a more
restrictive definition. Given a couple $(\cA, \cB)$ of an algebra and a
quotient manifold algebra, suppose there exists a Lie algebra $\fg$ and
an injective homomorphism of Lie algebras $\fg\rightarrow\der(\cA)$ such
that $\cB$ is the basic algebra for the operation of $\fg$ on $\cA$
(then $\fg\subset\hat{\fg}$). A connection on this triplet $(\cA,
\cB,\fg)$ is a spliting $\psi :\der(\cB) \rightarrow \fh$ of
$\cZ(\cB)$-modules compatible with the operation of $\fg$ in the sense
$[\fg, \psi(X)] = 0$ for all $X\in\der(\cB)$. Such a connection is also
given by a covariant projection $P: \fh \rightarrow \hat{\fg}$ where the
covariance means $[Y,P(X)] = P([Y,X])$ for all $Y\in\fg$ and $X\in\fh$.

In this situation, if $V$ is a linear space on which $\fg$ is
represented by $\eta: \fg \rightarrow \End(V)$, then one has an
associated central bimodule over $\cB$ defined by
\[ \cM_V = \{ a_i\otimes v^i \in \cA\otimes V \,\, /\,\, (Ya_i)\otimes
v^i + a_i\otimes \eta(Y)v^i = 0 \,\, \forall Y\in\fg \} \]
where the structure of bimodule over $\cB$ is localized on $\cA$.

\begin{proposition}
Let $\psi :\der(\cB) \rightarrow \fh$ be a connection on $(\cA,
\cB,\fg)$. Then, the mapping
\begin{eqnarray*}
\nabla_X^V : \cM_V &\rightarrow & \cM_V \\
a_i\otimes v^i &\mapsto & (\psi(X)a_i)\otimes v^i
\end{eqnarray*}
is well defined and is a connection on $\cM_V$. This is the associated
connection to $\psi$ on $\cM_V$.
\end{proposition}

\medskip
{\sl Proof:} $\nabla_X^V \cM_V \subset \cM_V$ because $[\fg,
\psi(X)]=0$. Other properties of $\nabla_X^V$ are immediate consequences
of the definition of $\psi$ as in
Proposition~\ref{connection}.~\hfill$\square$

\bigskip
In the case of a principal bundle, $\cM_V$ is the module over $\cB$ of
sections of the associated vector bundle for $(V,\eta)$. This module of
sections is considered here as the module of equivariant mapping
$P\rightarrow V$.

\bigskip
Let us now turn to a different problem. From the point of view of
characteristic classes (even if there is not yet such a theory for the
definition of connection used here), what is important in a connection
is its curvature. Given a couple $(\cA, \cB)$ of an algebra and a
quotient manifold algebra, suppose one has a central bimodule $\cM$ over
$\cA$ and a connection $\nabla$ on $\cM$, such that its curvature is
zero if one of its argument is in $\hat{\fg}$. Then one can transport
the connection on a central bimodule over $\cB$. Define the reduced
central bimodule over $\cB$
\[ \cM^{\hat{\fg}} = \{ m\in\cM \,\, /\,\, \nabla_X m = 0 \,\,
 \forall X\in\hat{\fg} \} \]

For any $\tilde{X}\in\der(\cB)$, take any $X\in\fh$ such that
$\rho(X)=\tilde{X}$. Then define, for any $m\in\cM^{\hat{\fg}}$,
\[ \widetilde{\nabla}_{\tilde{X}} m = \nabla_X m \]
Then, because the curvature of $\nabla$ is zero on $\hat{\fg}$, this is
a well defined mapping from $\cM^{\hat{\fg}}$ into itself. It is easy to
verify that $\widetilde{\nabla}$ is a connection, the curvature of which
is
\[ \tilde{R}(\tilde{X}, \tilde{Y})m = R(X,Y)m \]
for any $m\in\cM^{\hat{\fg}}$.

\bigskip
In the case where $\der(\cA)=\Int(\cA)$, it has been shown in
\cite{DVM2} that any central bimodule $\cM$ over $\cA$ admits the
canonical connection $\nabla_{ad(a)}m = am-ma$. The curvature of this
connection is zero.

Now, in the general case ($\der(\cA)\neq\Int(\cA)$), if one can take
this connection on $\Int(\cA)$ and a prolongement on $\der(\cA)$, then
the curvature is zero on $\Int(\cA)$. So one can hope to transport the
connection on a reduced module over $\cB=\cZ(\cA)$ while keeping the
same information on the curvature.

\section{Conclusion}

In this paper we have proposed definitions of the noncommutative
generalization of a submanifold and of a noncommutative quotient
manifold. Various examples and related constructions seems to give them
an importance for the study of derivations-based differential structures
of algebras. What must be notice is the different use of the two
generalizations of differential forms: $\underline{\Omega}_{\der}(\cA)$
and $\Omega_{\der}(\cA)$. This shows the importance to introduce various
generalizations of a commutative concept, adapted to different
situations.

\bigskip
\bigskip
{\bf Acknowledgements}

I would like to thank Michel Dubois-Violette for very helpful
discussions and John Madore for his kind interest.

\end{document}